\documentclass[11pt]{article}
\usepackage[latin1]{inputenc}
\usepackage[english]{babel}
\usepackage{geometry}
\usepackage{amsmath}
\usepackage{amssymb}
\usepackage{amsthm}
\usepackage{oldgerm}
\usepackage{amscd}
\usepackage{color}
\usepackage{hyperref}

\newtheorem{Theorem}{Theorem}[section]
\newtheorem{Lemma}{Lemma}[section]
\newtheorem{Cor}{Corollary}[section]
\theoremstyle{definition}
\newtheorem{Def}{Definition}[section]
\newtheorem{Ex}{Example}[section]
\newtheorem{Rem}{Remark}[section]

\newcommand{\Proof}{\noindent{\bfseries Proof : }} 
\newcommand{\Q}{\mathbb{Q}}
\newcommand{\N}{\mathbb{N}}
\newcommand{\Z}{\mathbb{Z}}
\newcommand{\IZ}{\textnormal{Int($\Z$)}}

\title{Factorization of integer-valued polynomials with square-free denominator}
\author{Giulio Peruginelli}
\date{\today}

\begin{document}
\maketitle
\begin{center}
\emph{{\small Dedicated to Marco Fontana on the occasion of his 65th birthday}}
\end{center}
\vskip0.3cm
\begin{abstract}
\noindent We  describe an algorithm  to compute the different factorizations of a given image primitive integer-valued polynomial $f(X)=g(X)/d\in\Q[X]$, where $g\in\Z[X]$ and $d\in\N$ is square-free, assuming that  the factorizations of $g(X)$ in $\Z[X]$ and $d$ in $\Z$ are known.
We translate this problem into a combinatorial one.
\end{abstract}


\vskip1cm
\section{Introduction}
It is well known that the ring of integer-valued polynomials $\IZ=\{f\in\Q[X]\,|\,f(\Z)\subset\Z\}$ is far from being a unique factorization domain. We know that the ring $\IZ$ is atomic (every non-zero non-unit of  $\IZ$ admits a factorization into irreducibles) and every non-zero non-unit has only has only finitely many factorizations into irreducibles (see \cite{SF}). In particular, this implies that $\IZ$ is a bounded factorization domain (the length of the different factorizations of a given element is bounded, see \cite[Prop. VI.3.2]{CaCh0}). Moreover, in \cite{CaChEl} it is shown that the ring has infinite elasticity, where the elasticity of a domain is defined as the supremum of the set of ratios between length of factorizations of non-zero non-units. We recall that the length of a factorization is the number of irreducible elements which appear in the factorization itself. Two factorizations into irreducibles of an element $x$ in a commutative ring $R$, say $x=r_1\cdot\ldots\cdot r_n$ and $x=s_1\cdot\ldots\cdot s_m$, are essentially the same if $n=m$ and after possibly re-indexing, $r_i$ is associated to $s_i$, for $i=1,\ldots,n$ (that is, there exists a unit $u_i\in R$ such that $r_i=u_i s_i$). Otherwise the two factorizations are essentially different (see [7]).

More recently, in \cite{SF} the following result is proved. Given a finite set  $S=\{n_1,\ldots,n_r\}$ of (non necessarily distinct) positive integers greater than 1, there exists an integer-valued polynomial $f(X)$ with   $r$ essentially distinct factorizations  into irreducibles of length $n_1,\ldots,n_r$, respectively. Hence, there are elements in the ring of integer-valued polynomials which admit distinct factorizations into irreducibles of arbitrary lengths. 

We propose here a new method to describe the essentially different factorizations into irreducibles of a given integer-valued polynomial, under the assumption that the denominator is square-free (we will treat the general case in a future work). We remark that in all the examples produced in \cite{SF} to exhibit polynomials with prescribed sets of lengths, only polynomials with square-free denominator appear (in \cite[Thm. 10]{SF} there is a polynomial with more than one prime in the denominator, in all the other results there is just one prime factor in the denominator). So, a treatment of this case has a certain interest. We begin by recalling some classical definitions.

\begin{Def}
The \textbf{content} of a polynomial $g(X)=\sum_{k=0,\ldots,n}a_kX^k\in\Z[X]$ is defined as the g.c.d. of its coefficients $a_k$. We denote the content of $g(X)$ by c$(g)$. A polynomial $g\in\Z[X]$ is called \textbf{primitive} if its content is equal to $1$. Given $f\in\IZ$, we denote by d$(f)$ the \textbf{fixed divisor} of $f$, that is the g.c.d. of the set of values $\{f(n)\,|\,n\in\Z\}$. An integer-valued polynomial $f(X)$ is said to be \textbf{image primitive} if its fixed divisor is equal to $1$. Let $p\in\Z$ be a prime. If $g\in\Z[X]$ we say that $g(X)$ is $p$-primitive if $p$ does not divide c$(g)$, that is, at least one of the coefficient of $g(X)$ is not divisible by $p$. If $f\in\IZ$, we say that $f(X)$ is $p$-image primitive if $p$ does not divide d$(f)$. 
\end{Def}

Given a polynomial $g\in\Z[X]$, the content of $g(X)$ is in general a proper divisor of the fixed divisor of $g(X)$: consider for example $g(X)=X(X-1)$ which is primitive but its fixed divisor is equal to $2$. Already in  \cite{ChMcCl} it is shown the important role played by the fixed divisor in the study of the factorizations of an integer-valued polynomial (see the results that we recall below). For example, the polynomial $g(X)=X^2+X+2$, which is irreducible in $\Z[X]$ (and consequently in $\Q[X]$ by Gauss Lemma), has fixed divisor equal to $2$, so that  in $\IZ$ we have the non trivial factorization $g(X)=2\cdot\frac{g(X)}{2}$.

We recall the following facts:
\begin{itemize}
\item[-] $\IZ$, $\Z[X]$ and $\Z$ share the same group of units, $\{\pm1\}$ (\cite[Lemma 1.1]{CaChEl}).
\item[-]  An irreducible integer $p$ stays irreducible in $\IZ$ (\cite[Lemma VI.3.1]{CaCh0}).
\item[-] $\IZ$ has no prime elements (\cite[Prop. 3.2]{ACaChS}).  
\item[-] If $g\in\Z[X]$ is $p$-primitive for some prime $p$ and $p$ divides the fixed divisor of $g(X)$, then $p\leq\deg(g)$ (this is due to Polya, see \cite[Thm. 3.1]{Nark} for a modern treatment). 
\item[-] Let $g\in\Z[X]$ be primitive. Then $g(X)$ is irreducible in $\IZ$ if and only if it is image primitive and irreducible in $\Z[X]$ (Chapman-McClain \cite[Thm. 2.6]{ChMcCl}). Hence,  an irreducible factor in $\IZ$ of an irreducible polynomial $g\in\Z[X]$ is  either a constant $c$ which is a divisor of $d(g)$ or $\frac{g(X)}{d(g)}$. 
\item[-] Given two integer-valued polynomials $f$ and $g$, we have $d(fg)\subset d(f)d(g)$ and in general we may not have an equality (see the above example $X(X-1)$). This is the main difference between fixed divisor and content: in fact, the content of the product of two polynomials $g_1(X)$ and $g_2(X)$ is equal to the product of the contents of $g_1(X)$ and $g_2(X)$ by Gauss Lemma. This is equivalent to the fact that a  primitive  polynomial $g\in \Z[X]$ is irreducible if and only if it is irreducible in $\Q[X]$; this sentence is no more true if we substitute the ring $\Q[X]$ with $\IZ$ (see the above example $g(X)=X^2+X+2$). By the above cited theorem of Chapman-McClain, we have to add the assumption that $g(X)$ is also image primitive. We notice that a factor of an image primitive polynomial is image primitive (\cite{ChMcCl}).
\end{itemize}
\vskip0.4cm
 
Given a polynomial $f\in\Q[X]$, we have $f(X)=g(X)/d$, for some uniquely determined $g\in\Z[X]$ and 
$d\in\N$  
such that $(d,{\rm c}(g))=1$ (we essentially use the fact that $\Z$ is UFD). For short, we call $d$ the denominator of $f(X)$ 
 and $g(X)$ the numerator of $f(X)$.  

We can further express $f(X)$ in the following way:
\begin{equation}\label{f(X)Z}
f(X)=\frac{g(X)}{d}=\frac{\prod_{i\in I}g_i(X)^{e_i}}{\prod_{k\in K}p_k^{e_k}}
\end{equation}
where $g(X)=\prod_{i\in I}g_i(X)^{e_i}$ is the unique irreducible factorization in $\Z[X]$ (the $g_i(X)$ may be possibly constant) and $d=\prod_{k\in K}p_k^{f_k}$ is the factorization of $d$ in $\Z$. Obviously, $f(X)$ is integer-valued if and only if $d$ divides the fixed divisor of $g(X)$, that is, for each $k=1,\ldots,m$, $p_k^{e_k}$ divides d$(g)$. Since $\IZ\subset\Q[X]$ and $\Q[X]$ is a UFD, any irreducible factor $h(X)$ of $f(X)$ in $\IZ$ is a rearrangement of the irreducible factors $g_i(X)$ of $g(X)$ and the prime factors $p_k$ of $d$, in such a way that we still have an integer-valued polynomial, that is: 
$$h(X)=\frac{g_1(X)}{d_1}=\frac{\prod_{i\in J}g_i(X)^{e_i'}}{\prod_{k\in T}p_k^{f_k'}}$$
where $J\subseteq I$, $T\subseteq K$, $e_i'\leq e_i,f_k'\leq f_k$ for each $i\in J$ and $k\in T$ and $h(X)$ is in $\IZ$, that is $d_1$ divides d$(g_1)$. It is already not clear how a general irreducible polynomial in $\IZ$ looks like (for polynomials $g\in\Z[X]$ which are irreducible in $\IZ$ see the above Theorem of Chapman-McClain). To our knowledge, the only characterization of such irreducibles is given by \cite[Cor. 2.9]{ChMcCl}, which largely relies on the problem of establishing the fixed divisor of a polynomial with integer coefficients. We will give a new characterization of the irreducible elements of $\IZ$ in the case of square-free denominator.

As we observed above, being image primitive is a necessary condition for an integer-valued polynomial to be irreducible.
 We will give a chacterization of image primitive polynomials.  
In particular, being image primitive implies that the numerator $g(X)$ is primitive (since we assume that the denominator $d$ is coprime with the content of $g(X)$). Notice that the converse of the previous statement does not hold, namely, if an integer-valued $f(X)$ is image primitive then it is not true in general that $f(X)$ is irreducible in $\IZ$. Consider for example $f(X)=X(X-1)^2/2$ which has the irreducible factor $X-1$ (by \cite[Cor. 2.2 \& Example 2.3]{CaChEl}, monic linear polynomials are irreducible in $\IZ$).

Let $p\in\Z$ be a fixed prime. We set 
$$I_p\doteqdot p\IZ\cap\Z[X]=\{g\in \Z[X] \mid p|d(g)\},$$
which is the ideal of polynomials in $\Z[X]$ whose fixed divisor is divisible by $p$. From \cite{Per} (but see also \cite[Chapt. 2, \textbf{18}, p. 22]{Dickson}) we know that
$$
I_p=(p,X^p-X)=\Big(p,\prod_{i=0,\ldots,p-1}(X-i)\Big)=\bigcap_{j=0}^{p-1}(p,X-j).
$$
The last intersection is precisely the primary decomposition of the ideal $I_p$ (see \cite[Lemma 2.2]{Per}). For $j=0,\ldots,p-1$, we set
$$\mathcal M_{p,j}\doteqdot(p,X-j)=\{g\in\Z[X] \mid p|g(j)\}.$$
The above intersection is actually equal to a product of ideals, since $\mathcal{M}_{p,j}$, for $j=0,\ldots,p-1$, are $p$ distinct maximal ideals in $\Z[X]$.  
More in general, if $n$ is a positive integer, we set
$$I_{p^n}\doteqdot p^n\IZ\cap\Z[X]$$
which is the ideal of polynomials whose fixed divisor is divisible by $p^n$. 
Clearly, a polynomial $f\in\Q[X]$ like in (\ref{f(X)Z}) is in $\IZ$ if and only if for every prime factor $p_k$ of the denominator $d$, the numerator $g(X)$ is in $I_{p_k^{e_k}}$.

In the next section we will introduce the notion of prime covering for 
 the set of irreducible factors of the numerator of an integer-valued polynomial $f(X)$. 
For each prime $p$ which appears in the denominator and for each irreducible polynomial $g\in\Z[X]$ which appears in the numerator, we look for the primary components of $I_p$ which contain $g(X)$. A subset of the irreducible factors $\{g_i(X)\}_{i\in I}$ of the numerator of $f(X)$ whose elements are contained in all the primary components of $I_p$ is called a $p$-covering. A $p$-covering is minimal if, whenever we remove an element, one of the primary components of $I_p$ is not covered by any of the polynomials left in the $p$-covering itself. 
 In the case of prime denominator, say $f(X)=\frac{g(X)}{p}=\frac{\prod_{i\in I}g_i(X)}{p}$, 
$f(X)$ is irreducible in $\IZ$ if and only if $\{g_i(X)\}_{i\in I}$ form a minimal $p$-covering. In the same way, 
if for a subset $J\subsetneq I$ we have $\prod_{i\in J}g_i(X)$ in $I_p$, then $f(X)$ is reducible in $\IZ$. If that choice is minimal in the above sense, then that factor is irreducible.

In the subsequent section, we generalize the previous results to the case of an integer-valued polynomial with square-free denominator. Finally, as an explicit example, we consider the case of an integer-valued polynomial with denominator equal to the product of two distinct primes.
\vskip1cm

\section{Integer-valued polynomials with prime denominator}\label{PC}

\subsection{Prime covering}

\begin{Def}\label{Cpg}
Let $g\in\Z[X]$ and $p\in\Z$ be a prime. We set
$$C_{p,g}\doteqdot\{j\in\{0,\ldots,p-1\} \mid p|g(j)\}.$$
\end{Def}
\vskip0.2cm
\noindent  Notice that the elements of $C_{p,g}$ correspond precisely to the primary components $\mathcal{M}_{p,j}=(p,X-j)$ of $I_p$ which contain $g(X)$. 
Observe that the set $C_{p,g}$ can be empty: for instance, take $g(X)=X^2+1$ and $p=3$. Obviously, $\#C_{p,g}\leq p$. Equivalently, we may consider the polynomial $\overline g\in(\Z/p\Z)[X]$ obtained by reducing the coefficients of $g$ modulo $p$. 
 A residue class $j\in \Z/p\Z$ is a root of $\overline g(X)$ if and only if the primary component $\mathcal{M}_{p,j}$ of $I_p$ contains $g(X)$.

We have the following result, which involves the family of sets $\{C_i\}_{i\in I}$ just defined. 
 We omit the proof, which follows directly from the definitions.
\begin{Lemma}\label{J}
Let $g(X)=\prod_{i\in I}g_i(X)$ be a product of polynomials in $\Z[X]$ and let $p$ be a prime. For each $i\in I$, let $C_i=C_{p,g_i}$. Then
$$g\in I_p\Leftrightarrow \bigcup_{i\in I}C_i=\{0,\ldots,p-1\}.$$
In particular, $g(X)$ is $p$-image primitive if and only if there exists $j\in\{0,\ldots,p-1\}$ such that no $C_i$ contains $j$.
\end{Lemma}
Notice that the condition $g(X)$ is $p$-image primitive is equivalent to $g\notin I_p$. 
Obviously we don't need to factor a given integer coefficient polynomial $g(X)$ in $\Z[X]$ in order to establish whether it is $p$-image primitive or not (just consider it modulo $p$ as we said above). By Polya's Theorem we cited in the introduction, it is sufficient to consider only those primes $p$ which are less or equal to the degree of $g(X)$. Anyway, for the study of the problem of the factorization in the ring $\IZ$ it is useful to write the statement as it is.

We give now the following definition.
\begin{Def}\label{pcovering}
Let $\mathcal G=\{g_i(X)\}_{i\in I}$ be a set of polynomials in $\Z[X]$. Let $p$ be a prime. For each $i\in I$ we set $C_i=C_{p,g_i}$. A \textbf{$p$-covering} for $\mathcal G$ (or just \textbf{prime covering}, if the prime $p$ is understood) is a subset $J$ of $I$ such that
$$\bigcup_{i\in J}C_i=\{0,\ldots,p-1\}.$$

We say that $J$ is \textbf{minimal} if no proper subset $J'$ of $J$ has the same property.
We will always assume that a given prime covering $J$ is \textsl{proper}, that is, for each $i\in J$ we have $C_i\not=\emptyset$.  

Two $p$-covering $J_1,J_2\subset I$ for $\mathcal{G}$ are disjoint if $J_1\cap J_2=\emptyset$.
\end{Def}
Notice that from a prime covering we can always extract a minimal prime covering, by discarding the redundant sets $C_i$. We may rephrase Lemma \ref{J} by saying that $f(X)=g(X)/p=\prod_{i\in I}g_i(X)/p$ belongs to $\IZ$ if and only if $I$ contains a $p$-covering for $\{g_i\}_{i\in I}$. 
A minimal $p$-covering can have $1$ element, for example consider the irreducible polynomial $X^p-X+p$. It has at most $p$ elements.
The problem to find such $p$-coverings has a combinatorial flavour.

The next example shows that given a minimal $p$-covering $J$, it does not follow that $\{C_i\}_{i\in J}$ forms a family of disjoint subsets of the residue classes modulo $p$. In fact, a polynomial $g_i(X)$ may belong to different primary components of $I_p$. If this is the case the degree of $g_i(X)$ has to be greater than one.

\vskip0.5cm
\begin{Ex}
\begin{equation}\label{notdisj}
f(X)=\frac{(X^2-X+3)(X^2+2)}{3}
\end{equation}
if we set $g_1(X)=X^2-X+3,g_2(X)=X^2+2$ we immediately see that
\begin{itemize}
\item $C_{3,g_1}=\{0,1\},C_{3,g_2}=\{1,2\}$
\item $C_{2,g_1}=\emptyset,C_{2,g_2}=\{0\}$
\end{itemize}
the second line implies that $2$ does not divide the fixed divisor of the numerator, that is $f(X)$ is $2$-image primitive (by Polya's Theorem, we check only those primes $p$ which are less or equal to the degree of $f(X)$).
We have that $C_{3,g_1}$ and $C_{3,g_2}$ covers the residue classes modulo $3$ and they have non trivial intersection. In particular, $I=\{1,2\}$ is a minimal $3$-covering. 
\end{Ex}

\vskip0.5cm
\subsection{Integer-valued polynomials which are $p$-image primitive}

We characterize now $p$-image primitive integer-valued polynomials, when the denominator is exactly divisible by a prime $p$ ( we denote this by $p\parallel d$ ).

Suppose that for a polynomial $f(X)$ as in (\ref{f(X)Z}) the denominator $d$ is equal to a prime $p$. If $f(X)$ is $p$-image primitive, then there exists $i\in I$ such that $e_i=1$, otherwise the fixed divisor of the numerator $g(X)$ is divisible by $p^n$, for some $n>1$. For example, $\frac{X(X-1)^2}{2}$ is $2$-image primitive, while $\frac{X^2(X-1)^2}{2}$ is not (the numerator has fixed divisor equal to $4$). However, this condition on the exponents of the irreducible factors $g_i(X)$ is not sufficient to ensure that $f(X)$ is $p$-image primitive, as the next example shows:
\begin{equation}\label{redmincov}
f(X)=\frac{(X^2+4)(X^2+3)}{2}
\end{equation}
The polynomial $f(X)$ is not $2$-image primitive since the numerator has fixed divisor equal to $4$ (modulo $2$, each factor at the numerator has a double root in $0$ and $1$, respectively).

Moreover, under the above assumption, the next lemma shows that all the minimal $p$-coverings must intersect in one spot. For $g(X)=\prod_{i\in I}g_i(X)\in\Z[X]$ and $J\subseteq I$, we set
$$g_J(X)\doteqdot\prod_{i\in J}g_i(X).$$
For each $i\in I$ we set $C_i=C_{p,g_i}$. By Lemma \ref{J}, for any subset $J\subseteq I$ we have $g_J\in I_p\Leftrightarrow J$ is a $p$-covering.

\begin{Lemma}\label{cond}
Let 
$$f(X)=\frac{\prod_{i\in I}g_i(X)}{d}$$
be in $\IZ$. Let $p$ be a prime factor of $d$ such that $p\parallel d$. Then $f(X)$ is $p$-image primitive if and only if the following condition holds: there exists a primary component $\mathcal{M}_{p,\overline j}$ of $I_p$, for some $\overline j\in\{0,\ldots,p-1\}$,  such that $g_{\overline i}\in\mathcal{M}_{p,\overline j}\setminus\mathcal{M}_{p,\overline j}^2$ for some $\overline i\in I$ and for all $i\in I$, $i\not=\overline i$, we have $g_{i}\notin\mathcal{M}_{p,\overline j}$.

If that condition holds, then for every minimal $p$-covering $J\subseteq I$, we have $\overline i\in J$. 

\end{Lemma}

\Proof Suppose $f(X)$ is $p$-image primitive. If for every $j\in\{0,\ldots,p-1\}$ there exist $i_1(j)\not=i_2(j)$ in $I$ such that $g_{i_1},g_{i_2}\in\mathcal{M}_{p,j}$, then we can form two disjoint $p$-coverings $J_t=\{i_t(j)\}_{j=0,\ldots,p-1}$, for $t=1,2$. By Lemma \ref{J} the polynomials $g_{J_1}$ and $g_{J_2}$ belong to $I_p$, thus their fixed divisor is divisible by $p$; since $g_I$ is divisible by $g_{J_1}\cdot g_{J_2}$, it has fixed divisor divisible by $p^2$, contradiction. So there exists $j'\in\{0,\ldots,p-1\}$ for which only one irreducible factor $g_{i'}(X)$ is in $\mathcal{M}_{p,j'}$. If $g_{i'}\notin\mathcal{M}_{p,j'}^2$ we are done. Suppose that is not the case. If for all the other $j$'s we have either more than one factor $g_i(X)$ in $\mathcal{M}_{p,j}$ or a factor $g_i(X)$ which belongs to $\mathcal{M}_{p,j}^2$ we get again to the same contradiction as before. Hence, there must be some  $\overline j\in\{0,\ldots,p-1\}$ for which the corresponding primary component $\mathcal{M}_{p,\overline j}$ of $I_p$ contains only one factor $g_{\overline i}(X)$. Moreover, $g_{\overline i}\notin\mathcal{M}_{p,\overline j}^2$.

Conversely, suppose there exists $\overline j\in\{0,\ldots,p-1\}$ as in the statement. If, for each $i\in I$, we set $C_i=C_{p,g_i}$ we have that $\overline{j}\notin C_i$ for all $i\not=\overline{i}$. Let $J\subseteq I$ be a minimal $p$-covering for $\{g_i\}_{i\in I}$ (we know that such a prime covering exists by Lemma \ref{J}).  Since by definition $\bigcup_{i\in J}C_i=\{0,\ldots,p-1\}$, and for all $i\in I$, $i\not=\overline i$, we have $C_i\not\ni\overline j$, it follows that $\overline i$ is contained in $J$. Notice that this proves the last statement of the Lemma. So there are no two disjoint $p$-coverings. Since $g_{\overline i}\notin\mathcal{M}_{p,\overline j}^2$ and $g_{\overline i}$ is the only factor of the numerator of $f(X)$ in $\mathcal{M}_{p,\overline j}$ we have that $g_J\notin I_{p^2}$. Since this holds for every minimal $p$-covering $J$, this concludes the proof of the lemma. $\Box$

\vskip0.5cm

\begin{Rem}
Under the assunptions of Lemma \ref{cond}, $f(X)$ is $p$-image primitive if and only if there exists a primary component $\mathcal{M}_{p,\overline j}$ of $I_p$ which contains one and only one irreducible factor $g_{\overline i}(X)$ of the numerator of $f(X)$ and $g_{\overline i}\notin \mathcal{M}_{p,\overline j}^2$. In particular, this means that only $C_{\overline i}$ contains $\overline j$. 

The last statement of Lemma \ref{cond} cannot be reversed, see for example (\ref{redmincov}). We have to add the hypothesis that for each minimal $p$-covering $J\subseteq I$ there exists at least one $i\in J$ such that $g_i\in \mathcal{M}_{p,j}\setminus \mathcal{M}_{p,j}^2$ for some $j\in J$. Equivalently, by the remarks after Definition \ref{Cpg}, we can say that for at least one residue classes $j$ modulo $p$, there is one and only one irreducible factor $g_i(X)$ which has a simple root modulo $p$ in $j$.

We can have more than one minimal $p$-covering, say $J_1,J_2\subseteq I$, provided they are not disjoint, as Lemma \ref{cond} says. For instance, consider the polynomial:
\begin{equation}\label{123}
f(X)=\frac{X(X-1)(X-2)}{2\cdot3}
\end{equation}
which is known to be irreducible (\cite[Example 2.8]{CaChEl}; in particular, $f(X)$ is image primitive). We set $g_{i+1}(X)=X-i$, for $i=0,1,2$. Then $J_1=\{1,2\}$ and $J_2=\{2,3\}$ are different minimal $2$-coverings, which are not disjoint.
\end{Rem}

\vskip0.3cm
\begin{Ex}
$$f(X)=\frac{X^2\cdot(X-1)\cdot(X^2+4)}{2}$$
in this example only $X-1$ belongs to $\mathcal{M}_{2,1}$ and moreover it does not belong to $\mathcal{M}_{2,1}^2$. Hence, the polynomial is $2$-image primitive.
\end{Ex}
\vskip0.3cm
\begin{Ex}
$$f(X)=\frac{X\cdot(X^2-2X+5)\cdot(X+6)}{2}$$
in this example, only $g(X)=X^2-2X+5$ belongs to $\mathcal{M}_{2,1}$. Moreover, $g\in\mathcal{M}_{2,1}^2$. No irreducible polynomial in the numerator belongs to $\mathcal{M}_{2,0}^2$, but there are two distinct factors, namely $X$ and $X+6$, which belong to $\mathcal{M}_{2,0}$. Hence, $f(X)$ is not $2$-image primitive, since the fixed divisor of the numerator is $4$. So it is not sufficient to have a unique $\overline i\in I$ such that $g_{\overline i}\in\mathcal{M}_{p,\overline j}$. We must also take care of the exact power of the maximal ideal $\mathcal{M}_{p,j}$ to which each polynomial $g_{i}(X)$ belongs to.
\end{Ex}
\vskip0.8cm

\subsection{Irreducible integer-valued polynomials}

Suppose that an integer-valued polynomial $f(X)$ is of the form
\begin{equation}\label{f(X)p}
f(X)=\frac{g(X)}{p}=\frac{\prod_{i\in I}g_i(X)}{p}
\end{equation}
where, for $i\in I$, $g_i\in\Z[X]$ is irreducible. The fact that $f\in\IZ$ is image primitive amounts to saying that $d(g)$ is equal to $p$.
Since $f\in\IZ$, by Lemma \ref{J} there exists a $p$-covering $J\subseteq I$ for $\{g_i(X)\}_{i\in I}$. The next lemma establishes that $f(X)$ is irreducible in $\IZ$ if and only if $I$ is a minimal $p$-covering. 
\begin{Lemma}\label{minpcov}
An image primitive polynomial $f(X)=g(X)/p$ in $\IZ$ as in (\ref{f(X)p}) is irreducible in $\IZ$ if and only if there is no proper subset $J$ of $I$ such that $\bigcup_{j\in J}C_j=\{0,\ldots,p-1\}$ (that is, $I$ is a minimal $p$-covering).
\end{Lemma}
\vskip0.3cm
\Proof Suppose there exists $J\subsetneq I$ such that $J$ is a $p$-covering. Then
$$f(X)=\frac{g_J(X)}{p}\cdot g_{I\setminus J}(X)$$
is a non-trivial factorization of $f(X)$ in $\IZ$, because the first factor is integer-valued by Lemma \ref{J} and the second one is in $\Z[X]\subset\IZ$. 

Conversely, if $f(X)$ is reducible in $\IZ$, then there exist non-constant $g,h\in\IZ$ such that $f(X)=h_1(X)h_2(X)$ (because we are assuming $f(X)$ to be image primitive). Since $p$ must appear in the denominator of one of the two factors, say $h_1(X)$, then for some $\emptyset\not=J\subsetneq I$ we have $h_1(X)=g_J(X)/p$ and consequently $h_2=g_{I\setminus J}\in\Z[X]$. Since $h_1\in\IZ$, by Lemma \ref{J} $J$ is a $p$-covering (notice that $h_1\in\IZ\Leftrightarrow g_J\in I_p$). $\Box$
\vskip0.3cm 
Notice that Lemma \ref{minpcov} does not hold without assuming $f(X)$ to be image primitive, as example (\ref{redmincov}) shows. By the arguments we have just given, we deduce that every factorization of an image primitive integer-valued polynomial with prime denominator $f(X)=\frac{g(X)}{p}$ is of the form $f(X)=\frac{g_J(X)}{p}\cdot g_{I\setminus J}(X),$ for some $J\subseteq I$ minimal $p$-covering.
Notice that the number of irreducible factors of the previous factorization in $\IZ$ is $1+\#(I\setminus J)$. The assumption that $f(X)$ is image primitive implies that for each such a minimal $p$-covering $J$, the set $I\setminus J$ does not contain a $p$-covering.
\vskip1cm

\section{Integer-valued polynomials with square-free denominator}

The main problem in the general case of more than one prime factor in the denominator $d$ of an integer-valued polynomial $f(X)$ is that each irreducible factor $g_i(X)$ of the numerator of $f(X)$ may belong to different primary components $\mathcal{M}_{p_k,j}$ of $I_{p_k}$, where $\{p_k\}_{k\in K}$ are the different prime factors of $d$.

As already remarked in \cite{Lewis}, this phenomenon has the effect that if $p$,$q$ are two distinct primes, then it does not follow that $I_p\cdot I_q=I_{pq}$: for example, $g(X)=X(X-1)(X-2)$ belongs to $I_{2\cdot 3}$ (see (\ref{123})), but it cannot be expressed as a product of a polynomial in $I_2$ and a polynomial in $I_3$. This is due to the fact that the only minimal $3$-covering $J=\{1,2,3\}$ is equal to the set $I$ itself, so in particular it has non-zero intersection with any possible $2$-covering (we saw that there are only two of them). Hence, in the next subsection, we are lead to give this globalizing definition.
\vskip0.3cm
\subsection{Family of minimal $\mathcal{P}$-coverings}
\begin{Def}
Let $\mathcal G=\{g_i(X)\}_{i\in I}$ be a set of polynomials in $\Z[X]$ and let $\mathcal P=\{p_k\}_{k\in K}$ be a set of distinct prime integers. A \textbf{family of minimal $\mathcal P$-coverings for $\mathcal G$} is a family of sets $\{J_k\}_{k\in K}$ such that for each $k\in K$, $J_k\subseteq I$ is a minimal $p_k$-covering for $\mathcal G$.
\end{Def}
\vskip0.3cm
Let $f\in\Q[X]$ be as in (\ref{f(X)Z}). If $f(X)$ is an integer-valued polynomial, then by Lemma \ref{J} there exists a family of minimal $\mathcal P=\{p_k\}_{k\in K}$-coverings for $\mathcal G=\{g_i(X)\}_{i\in I}$.

We can now formulate a proposition, which gives a criterion for an integer-valued polynomial to be irreducible, in the case that the denominator is square-free. This is a first step to determine explicitly all the factorizations of a given element in the ring $\IZ$.

Firstly we set some notations. Let 
\begin{equation}\label{gIpK}
f(X)=\frac{g(X)}{d}=\frac{\prod_{i\in I}g_i(X)}{\prod_{k\in K}p_k}
\end{equation}
be a polynomial in $\Q[X]$, with $p_k$ distinct prime integers, $g_i\in\Z[X]$ irreducible polynomials.  
Notice that the condition that $f(X)$ is integer-valued is equivalent to $g_I(X)=\prod_{i\in I}g_i(X)\in\bigcap_{k\in K}I_{p_k}$.
We set $\mathcal G=\{g_i(X)\}_{i\in I}$ and $\mathcal P=\{p_k\}_{k\in K}$. As in the previous section, given $J\subseteq I$ we set $g_J(X)\doteqdot\prod_{i\in J}g_i(X)$. Notice that if $J_1\subseteq J_2\subseteq I$ we have that $g_{J_1}(X)$ divides $g_{J_2}(X)$ in $\Z[X]$ (and so in $\IZ$). Similarly, for a subset $T\subseteq K$ we set 
$$d_T\doteqdot\prod_{k\in T}p_k$$
($d_K=d$). With these notations, a factor of $f(X)$ is of the form:
$$h(X)=\frac{g_J(X)}{d_T}$$
for some $J\subseteq I$ and $T\subseteq K$.

Finally, if $T\subseteq K$ and $\mathcal J=\{J_k\}_{k\in K}$ is a family of minimal $\mathcal P$-coverings for $\mathcal G$, we set 
$$I_{\mathcal J,T}\doteqdot\bigcup_{k\in T} J_{k}.$$
Notice that, if $T_1,T_2\subseteq K$ are two disjoint subsets, then $I_{\mathcal J,T_1\cup T_2}=I_{\mathcal J,T_1}\cup I_{\mathcal J,T_2}$.

\subsection{Irreducible integer-valued polynomials}

\vskip0.3cm
\begin{Theorem}\label{irr}

Let
$$f(X)=\frac{g(X)}{d}=\frac{\prod_{i\in I}g_i(X)}{\prod_{k\in K}p_k}$$
be an image primitive integer-valued polynomial.  Let $\mathcal P=\{g_i(X)\}_{i\in I}$ and $\mathcal G=\{p_k\}_{k\in K}$. We suppose that the polynomials $g_i(X)$ are irreducible in $\Z[X]$ and that the $p_k$ are distinct prime integers. Then $f(X)$ is irreducible in $\IZ$ if and only if the following holds: for every family $\mathcal J=\{J_k\}_{k\in K}$ of minimal $\mathcal P$-coverings for $\mathcal G$ we have
 \begin{itemize}
  \item[i)]  $I=I_{\mathcal J,K}$.
  \item[ii)] there is no non-trivial partition $K=K_1\dot\cup K_2$ such that $I_{\mathcal J,K_1}\cap I_{\mathcal J,K_2}=\emptyset$.
 \end{itemize}

\end{Theorem}

Notice that condition i) implies that for each $i\in I$ there exists $k\in K$ such that $C_{p_k,g_i}\not=\emptyset$, so that each of the $g_i$'s belongs to at least one of the primary components $\mathcal{M}_{p_k,j}$ of some of the ideals $I_{p_k}$. Moreover,  condition ii) says that the union of the elements $J_k$ of the family $\mathcal J$ cannot be partitioned (in a sense we will make precise soon).
We will treat the case $\mathcal P=\{p_1,p_2\}$ as an example in section \ref{2primes}.

\vskip0.5cm

\Proof Suppose $f\in\IZ$ irreducible. Let $\mathcal J=\{J_k\}_{k\in K}$ be a family of minimal $\mathcal P$-coverings for $\mathcal G$ (it exists because of Lemma \ref{J}). If $I$ strictly contains $I_{\mathcal J,K}$ then there exists $t\in I$ which is not contained in any $J_k$ (equivalently, $J_k\subseteq I\setminus\{t\}$ for every $k\in K$). This means that $g_t(X)$ divides $f(X)$ in $\IZ$, because we have
$$f(X)=g_t(X)\cdot\frac{g_{I\setminus\{t\}}(X)}{d}$$
and the second factor is integer-valued, since for each $k\in K$ we have $g_{J_k}(X)\in I_{p_k}$ (see Lemma \ref{J}). Hence, for all such $k$'s,  we have $g_{I\setminus\{t\}}(X)\in I_{p_k}$, since $J_k\subseteq I\setminus\{t\}$. This is a contradiction, hence condition i) holds.

If we have a non-trivial partition $K=K_1\dot\cup K_2$ such that $I_1\doteqdot I_{\mathcal J,K_1}$ and $I_2\doteqdot I_{\mathcal J,K_2}$ are disjoint, then
$$f(X)=\frac{g_{I_1}(X)}{d_{K_1}}\cdot\frac{g_{I_2}(X)}{d_{K_2}}.$$
Notice that for every $k_1\in K_1$ we have $g_{J_{k_1}}(X)\in I_{p_{k_1}}$ (again by Lemma \ref{J}) and $g_{J_{k_1}}(X)$ divides $g_{I_1}(X)$ in $\Z[X]$, since $J_{k_1}\subset I_1$. This implies that $g_{I_1}(X)/d_{K_1}$ is integer-valued. Similarly, the second factor is integer-valued, too. That would be a non-trivial factorization of $f(X)$, which is a contradiction.

Conversely, suppose that for every family $\mathcal J=\{J_k\}_{k\in K}$ of minimal $\mathcal P$-coverings for $\mathcal G$ conditions i) and ii) hold. Since $f(X)$ is image primitive, there is no non-unit in $\Z$ which divides $f(X)$ in $\IZ$. If $f(X)$ is reducible in $\IZ$ we have $f(X)=h_1(X)h_2(X)$, where $h_1,h_2\in\IZ$ are not constant. Since $\IZ\subset\Q[X]$ we have
$$h_i(X)=\frac{g_{I_i}(X)}{d_{K_i}}$$
for some $I_i\subseteq I$ and $K_i\subseteq K$, for $i=1,2$. Necessarily, $I_1,I_2$ are disjoint and $I_1\cup I_2=I$. Similarly, $K_1$ and $K_2$ are disjoint and $K_1\cup K_2=K$. Suppose that one of the $K_i$, say $K_2$, is empty. Then, by Lemma \ref{J} for each $k\in K_1=K$ there exists a minimal $p_k$-covering $J_k\subseteq I_1$. We set $\mathcal{J}=\{J_k\}_{k\in K}$. By definition, the family $\mathcal{J}$ is a minimal $\mathcal P$-coverings for $\mathcal G$. In particular, $I_{\mathcal J,K}\subseteq I_1$, because each of the $J_k$'s is a subset of $I_1$. Because of i) we have that $I=I_{\mathcal J,K}$, so that $I=I_1$ and consequently $I_2=\emptyset$, since $I_1$ and $I_2$ are disjoint. This means that $h_2(X)$ is a unit.

Suppose now that $K_i\not=\emptyset$, for $i=1,2$. This fact also leads us to a contradiction. In fact, by Lemma \ref{J}, for each $i=1,2$ and for each $k_i\in K_i$ there exists a minimal $p_{k_i}$-covering $J_{k_i}\subseteq I_i$. We set $\mathcal J=\{J_k\}_{k\in K}$, which is a family of minimal $\mathcal P$-coverings for $\mathcal G$.
In particular, $I_{\mathcal J,K_i}\subseteq I_i$. By condition i) on $\mathcal J$ we have that 
$$I=I_{\mathcal J,K}=I_{\mathcal J,K_1}\dot\cup I_{\mathcal J,K_2}.$$
Since $I_1\cup I_2=I$, we get $I_{\mathcal J,K_i}=I_i$ for $i=1,2$, which is in contradiction with condition ii). $\Box$
\vskip0.5cm

\begin{Ex}\label{Esempio}

It is not sufficient that conditions i) and ii) of Theorem \ref{irr} hold only for one family $\{J_k\}_{k\in K}$ of minimal $\mathcal{P}$-coverings. For instance, let us consider
$$f(X)=\frac{(X-1)\cdot(X-2)\cdot(X-3)\cdot(X-9)}{2\cdot 3}$$
then if $g_i(X)=X-i$, for $i=1,2,3$, $g_4(X)=X-9$ and $I=\{1,2,3,4\}$, we have that
\begin{itemize}
\item[-] $J_2=\{2,1\}$, $J_2'=\{2,3\}$ and $J_2''=\{2,4\}$ are the minimal $2$-coverings.
\item[-] $J_3=\{1,2,3\}$ and $J_3'=\{1,2,4\}$ are the minimal $3$-coverings.
\end{itemize}
We have that $\mathcal J=\{J_2'',J_3\}$ is a family of minimal $\mathcal P$-coverings for $\mathcal G$ which satisfies both conditions i) and ii) but the polynomial is not irreducible, since $X-9$ divides $f(X)$ in $\IZ$. In fact, the family $\mathcal{J}'=\{J_2',J_3\}$ of $\mathcal P$-coverings for $\mathcal G$ does not satisfy condition i) of the proposition.
\end{Ex}

\begin{Rem}
From Theorem \ref{irr} we see that each family of minimal $\mathcal P$-coverings for $\mathcal G$ determines a (possibly trivial, like for $\mathcal J$ in Example \ref{Esempio}) factorization for $f(X)$ in $\IZ$. Conversely, every non-trivial factorization determines a family of minimal $\mathcal P$-coverings for $\mathcal G$ which can be partitioned in the following sense:

\begin{Def}\label{family partition}
We say that a family $\mathcal{J}$ of minimal $\mathcal{P}$-coverings for $\mathcal{G}$ is \textbf{partitionable} if there exist a partition for $K$, say $K=\dot\bigcup_{j\in\mathcal{I}}K_j$ such that the sets $\{I_{\mathcal{J},K_j}=\bigcup_{k\in K_j}J_k \mid j\in\mathcal{I}\}$ are disjoint.
\end{Def}

However, notice that different families of minimal $\mathcal{P}$-coverings may give the same factorization for $f(X)$. For instance, in the Example \ref{Esempio}, there are six possible such families (we have to pair each minimal $2$-covering with a minimal $3$-covering). The family $\mathcal{J}''=\{J_2,J_3\}$ gives the same factorization as $\mathcal{J}'$. This depends on the fact that $I_{\mathcal{J}',K}$ and $I_{\mathcal{J}'',K}$ are equal.

\end{Rem}
\vskip0.5cm
\begin{Cor}
Let $f(X)$ be as in the assumptions of Theorem \ref{irr}. If there exists $\overline{k}\in K$ such that $I$ is a minimal $p_{\overline{k}}$-covering, then $f(X)$ is irreducible in $\IZ$.
\end{Cor}
\Proof We retain the notations of Theorem \ref{irr}. Let $\mathcal{J}$ be a family of minimal $\mathcal{P}$-coverings for $\mathcal{G}$. Notice that $I$ is the only minimal $p_{\overline{k}}$-covering, so $I\in \mathcal{J}$ and consequently $I=I_{\mathcal{J},K}$ and $\mathcal{J}$ is not a partitionable family. Hence, the conditions i) and ii) of Theorem \ref{irr} are satisfied for every family of minimal $\mathcal{P}$-coverings, so $f(X)$ is irreducible in $\IZ$. $\Box$
\vskip0.5cm
In particular, this corollary shows again that the polynomial in (\ref{123}) is irreducible. 
The condition of the previous corollary is not necessary, see (\ref{2nd}) below for an example.

\subsection{The algorithm of factorization in $\IZ$}

The next corollary shows explicitly how to obtain a non-trivial factorization of an integer-valued polynomial $f(X)$ as in (\ref{gIpK}) from a partitionable family of minimal $\mathcal{P}$-coverings for $\mathcal{G}$. We know from the proof of Theorem \ref{irr} that every such factorization is obtained in this way.

We recall that we are assuming $f(X)$ to be image primitive and the denominator of $f(X)$ to be square-free. 

Schematically we are doing the following steps:
\begin{itemize}
\item[i)] For each $k\in K$ and for each $i\in I$ we determine the sets $C_{p_k,g_i}$.
\item[ii)] Afterwards for each $k\in K$ we find all the minimal $p_k$-coverings $J_k$, by grouping together the sets $C_{p_k,g_i}$.
\item[iii)] Then for each $k\in K$ we choose one of the minimal $p_k$-coverings we found at point ii) and we define the family $\mathcal J=\{J_k\}_{k\in K}$ of minimal $\mathcal P$-coverings for $\mathcal G$.
\end{itemize}  
\vskip1cm
\begin{Cor}
 Let 
 $$f(X)=\frac{\prod_{i\in I}g_i(X)}{\prod_{k\in K}p_k}=\frac{g_I(X)}{d_K}$$
 be an image primitive, integer-valued polynomial, where $p_k$ are distinct prime integers, $g_i\in\Z[X]$ distinct and irreducible.
 
 Every factorization of $f(X)$ in $\IZ$ is obtained in the following way:
 
 let $\mathcal J=\{J_k\}_{k\in K}$ be a family of minimal $\mathcal P$-coverings for $\mathcal G$ which is partitionable, say $K=\dot\bigcup_{j\in \mathcal{I}} K_j$, so that the sets $I_j\doteqdot I_{\mathcal J,K_j}=\bigcup_{k\in K_j}J_k$, for $j\in\mathcal{I}$, are disjoint and for each $j\in\mathcal{I}$ the integer-valued polynomial $g_{I_j}(X)/d_{K_j}$ satisfies the conditions of Theorem \ref{irr} (so that each of them is irreducible). We set $I'\doteqdot\bigcup_{j\in \mathcal{I}}I_j$. Then
 $$f(X)=g_{I\setminus I'}(X)\cdot\prod_{j\in \mathcal{I}}\frac{g_{I_j}(X)}{d_{K_j}}$$
 is a factorization of $f(X)$ in $\IZ$ and every one of them is obtained in that way. Notice that in the previous factorization we have $\#(I\setminus I')+\#\mathcal{I}$ irreducible factors.
\end{Cor}

\vskip1cm
\subsection{Case $d=p_1\cdot p_2$}\label{2primes}

Let $p_1,p_2\in\Z$ be distinct primes. We consider an image primitive integer-valued polynomial of the following form:
\begin{equation}\label{f2}
f(X)=\frac{g(X)}{p_1 p_2}=\frac{\prod_{i\in I}g_i(X)}{p_1 p_2}
\end{equation}
This amounts to saying that the fixed divisor $d(g)$ is equal to $p_1 p_2$. By Lemma \ref{J}, for $k=1,2$, there exists a $p_k$-covering $J_k$ for $\{g_i(X)\}_{i\in I}$. For each $i\in I$ and for each $k=1,2$ we consider the sets $C_{p_k,g_i}$ as defined in section \ref{PC}. 
With the notation we introduced so far, we can have two different kind of factorization of $f(X)$. One possible factorization is:
\begin{equation}\label{1}
f(X)=\frac{g_{J_1}(X)}{p_1}\cdot\frac{g_{J_2}(X)}{p_2}\cdot\prod_{i\in I\setminus J_1\cup J_2}g_i(X)
\end{equation}
for some $J_{1},J_{2}\subseteq I$, where, for $k=1,2$, $g_{J_k}(X)/p_k\in\IZ$ is irreducible. By Lemma \ref{minpcov}, this corresponds to the fact that, for $k=1,2$, $J_k$ is a minimal $p_k$-covering. Obviously, $J_1$ and $J_2$ are disjoint. 

Another possible factorization is
\begin{equation}\label{2}
f(X)=\frac{g_{J}(X)}{p_1p_2}\cdot\prod_{i\in I\setminus J}g_i(X)
\end{equation}
for some $J\subseteq I$. In this factorization $g_{J}(X)/(p_1p_2)\in\IZ$ is irreducible. 

By Lemma \ref{J}, since $g_{J}(X)/(p_1p_2)$ is integer-valued then for each $k=1,2$, $J$ contains a minimal $p_k$-covering $J_k$. By Theorem \ref{irr}, the fact that $g_{J}(X)/(p_1p_2)$ is irreducible in $\IZ$ is equivalent to saying that $J=J_1\cup J_2$ (otherwise, we can factor out some $g_i(X)$ from it) and $J_1\cap J_2\not=\emptyset$ (otherwise we fall in the previous case (\ref{1})). It is not true that for some $k=1,2$ we must have $I=J_k$, like example (\ref{2nd}) below shows.

In \cite[Example 3.6]{ChMcCl} the authors construct an integer-valued polynomial which has two distinct factorizations as in (\ref{1}) and (\ref{2}). Now we give other two explicit examples: in the first one only the factorization as in (\ref{1}) occurs, in the second one we give an irreducible polynomial in $\IZ$ of the form $g(X)/(p_1p_2)$.
\vskip0.5cm
\begin{Ex}
\begin{align*}
f(X)&=\frac{(X^2+12)(X^2+2)(X^2+10)(X^2+16)(X^2+4)}{3\cdot5}\\
    &=\frac{(X^2+12)(X^2+2)}{3}\cdot\frac{(X^2+10)(X^2+16)(X^2+4)}{5}
\end{align*}
the second line is the only factorization in $\IZ$ that $f(X)$ can have, since if we put
$g_1(X)=X^2+12,g_2(X)=X^2+2,g_3(X)=X^2+10,g_4(X)=X^2+16,g_5(X)=X^2+4$ we have:
$$\begin{array}{lllll}
C_{3,g_1}=\{0\},&C_{3,g_2}=\{1,2\},&C_{3,g_3}=\emptyset,&C_{3,g_4}=\emptyset,&C_{3,g_5}=\emptyset\\
C_{5,g_1}=\emptyset,&C_{5,g_2}=\emptyset,&C_{5,g_3}=\{0\},&C_{5,g_4}=\{2,3\},&C_{5,g_5}=\{1,4\}
\end{array}$$

so in $I=\{1,\ldots,5\}$ we only have one $3$-covering $J_3=\{1,2\}$ and only one $5$-covering $J_5=\{3,4,5\}$, and they are disjoint. It is easy to check that $2$ and $7$ do not divide the fixed divisor of the numerator of $f(X)$.
\end{Ex}

\begin{Ex}
\begin{equation}\label{2nd}
f(X)=\frac{X(X^2+2)(X^2+16)(X^2+4)}{3\cdot5}
\end{equation}
so if $g_1(X)=X,g_2(X)=X^2+2,g_3(X)=X^2+16,g_4(X)=X^2+4$ we have:
$$\begin{array}{llll}
C_{3,g_1}=\{0\},&C_{3,g_2}=\{1,2\},&C_{3,g_3}=\emptyset,&C_{3,g_4}=\emptyset\\
C_{5,g_1}=\{0\},&C_{5,g_2}=\emptyset,&C_{5,g_3}=\{2,3\},&C_{5,g_4}=\{1,4\}
\end{array}$$
Then by Theorem \ref{irr} $f(X)$ is irreducible in $\IZ$ since $J_3=\{1,2\}$ is the only minimal $3$-covering, $J_5=\{1,3,4\}$ is the only minimal $5$-covering, $I=J_3\cup J_5$ and $J_3\cap J_5\not=\emptyset$. Notice that $J_3\subsetneq I, J_5\subsetneq I$. It is easy to check that $2$ and $7$ do not divide the fixed divisor of the numerator $g(X)$ of $f(X)$. In particular, $f(X)$ is image primitive, that is $d(g)=3\cdot5$.
\end{Ex}

\begin{Ex}
As another application of Theorem \ref{irr} we consider the polynomial:
$$f(X)=\frac{X\cdot(X^2+1)\cdot(X^2+X+1)\cdot(X^2+2X+4)}{2\cdot3}$$
and let $g_1(X)=X$, $g_2(X)=X^2+1$, $g_3(X)=X^2+X+1$, $g_4(X)=X^2+2X+4$. Then
\begin{itemize}
 \item[-] $J_2=\{1,2\}$ and $J_2'=\{2,4\}$ are the minimal $2$-coverings.
 \item[-] $J_3=\{1,3,4\}$ is the only minimal $3$-covering.
\end{itemize}
So $\mathcal J=\{J_2,J_3\}$ is a family of minimal $\mathcal P$-coverings of $\mathcal G$ such that $J_2\subsetneq I$, $J_3\subsetneq I$. The same holds for $\mathcal J'=\{J_2',J_3\}$. The polynomial is irreducible by Theorem \ref{irr}: if we consider $\mathcal J$ we have $I=J_2\cup J_3$ and $J_2\cap J_3\not=\emptyset$. The same holds for $\mathcal J'$.
\end{Ex}
\vskip0.2cm
Our method can be easily generalized to the case of denominator divisible by prime powers $p^n$ such that $n\leq p$, since in this case, by \cite[Proposition 3.1]{Per}, the primary components of the ideal $I_{p^n}$ are just the $n$-th power of the maximal ideals $\mathcal{M}_{p,j}$, for $j=0,\ldots,p-1$. In general, a further study of the primary components of the ideal $I_{p^n}$ is needed.

\vskip0.6cm
\subsection*{\textbf{Acknowledgments}}
I would like to thank the referee for his/her valuable suggestions. The author was supported by the Austrian Science Foundation (FWF), Project Number P23245-N18 and also by INDAM - GNSAGA.
\vskip0.3cm

\end{document}